\documentclass
       {article}
\usepackage{theorem}
\usepackage{oldgerm}
\usepackage{xspace}
\usepackage{enumerate}
\usepackage{latexsym}
\usepackage{amssymb}
\usepackage{amsmath}
\usepackage{epsfig}

\theoremstyle{plain}
        \newtheorem{thm}{Theorem}[section]
        \newtheorem{lem}[thm]{Lemma}

        \newtheorem{prop}[thm]{Proposition}
        \newtheorem{cor}[thm]{Corollary}
        
        \newtheorem{probl}[thm]{Problem}
        {\theorembodyfont{\rmfamily}
                \newtheorem{defn}[thm]{Definition}
                \newtheorem{rem}[thm]{Remark}
                \newtheorem{exa}[thm]{Example}
        }

\def\negyzet{\vbox{\hrule
                   \hbox{\vrule\kern2pt
                         \vbox{\kern2pt\kern2pt
                        }\kern2pt\vrule
                  }\hrule
                   }
             }

\renewcommand{\rho}{\varrho}

\newcommand{\proof}{\par\noindent {\sl Proof: \/}}

\newcommand{\proofof}[1]{\par\noindent {\sl Proof of #1: \/}}
\newcommand{\proofend}{~\rule{2mm}{3mm}}
\newcommand{\rnh}{$\rho$--neighborhood\xspace}
\newcommand{\rzn}{$\le\rho$--zone\xspace}
\newcommand{\Er}{E_{\le\rho}}
\newcommand{\utacfC}{up to a constant factor $C$\xspace}
\newcommand{\krl}{$k$--rich line\xspace}

\newcommand{\EP}{Erd\H{o}s\xspace}

\newcommand{\defeq}{\buildrel{\scriptstyle\rm def}\over=}

\def\dist{\hbox{\text{dist}}}
\def\segm#1{\overline{#1}}

\newcommand{\ds}{\textstyle}

\def\SzTr{Szemer\'edi--Trotter\xspace}

\def\pni{\par\noindent}
\newcommand{\p}{^{\prime}}

\newcommand{\A}{{\mathcal A}}

\newcommand{\Cell}{{\mathcal C}}

\newcommand{\J}{{\mathcal J}}
\renewcommand{\L}{{\mathcal L}}

\renewcommand{\P}{{\mathcal P}}

\newcommand{\R}{\mathbb R}

\begin{document}

\date{}

\title {On the Dimension of Finite Point Sets I.\\
       An Improved Incidence Bound for Proper 3D sets.}
\author{Gy\"orgy Elekes}
\maketitle

\begin{abstract}\noindent
We improve the well-known  \SzTr incidence bound
for proper 3--dimensional point sets (defined appropriately).
\end{abstract}

\begin{section}{Introduction}

\begin{subsection}{The \SzTr incidence bound in the plane}

The following estimate was conjectured by \EP and proven by \SzTr
for incidences of points and lines.
\begin{prop}[\SzTr Theorem]
\label{SzTrIncidenceProp}
\ The maximum number \ $I(n,m)$ of  incidences between $n$ points and $m$
straight lines in  the Euclidean plane satisfies
\[
I(n,m)=O(n^{2/3}m^{2/3} + n + m).
\]
As a special case, given a set $\P\subset\R^2$ of $n$ points, 
the number of $k$--rich lines (which contain at least $k$ points 
of $\P$) is bounded by
\[
C\cdot\max\Bigl\{\frac{n^2}{k^3},\frac{n}{k}\Bigr\},
\]
for an absolute constant $C>0$ and any $2\le k\le n$.
\end{prop}
\proof
see \cite{SzTr:83} and,  for a simple proof, \cite{Leslie:SzT}.
\pni
It is also true that these bounds are sharp, apart from the constant factors.
\smallskip\par
How much does the situation change if we consider point sets (and straight 
lines) in higher dimensional spaces?
On the one hand, the foregoing bounds still apply, as shown by a
projection to a generic plane.
On the other hand, no better bound can be stated, since 
any planar point set which attains the order of
magnitude in the \SzTr bounds, can be considered as 
a subset of $\R^3$ (or that of $\R^d$).
However, one might have the feeling that the real question
would be to consider \emph{proper } 3--dimensional
(or $d$--dimensional) sets.
\par
The main goal of this paper is to improve the \SzTr bound(s)
for proper 3--dimensional point sets (defined appropriately
in the next section).
\end{subsection}

\begin{subsection}{Proper $d$--dimensional point sets}
\label{ProperDimDefSection}

Let $H$ be a (finite) set of $n$ planes in $\R^3$
(or, in general, of hyperplanes in $\R^d$).
They cut the space into at most 
${n \choose 3} + {n \choose 2} + {n \choose 1} + {n \choose 0} \sim n^3$
open convex cells
(and into $\le \sum_{i=0}^d {n \choose i} \sim n^d$ in $\R^d$),
with equality iff $H$ is in general position,
i.e., if any three (in general, any $d$) share exactly one common point.
\begin{defn}
\label{ProperDimDef}
A set of $N$ points is \emph{proper $d$--dimensional 
\utacfC,\/} (for short, ``proper $d$--D'') 
if it  can be  cut into singletons by
at most $ C\root{d}\of{N}$ appropriate hyperplanes.
\end{defn}
\end{subsection}

\begin{subsection}{The main result}

\begin{thm}
\label{MainThm}
Assume that a set $\P\subset\R^3$ of $N$ points is proper 3--dimensional,
\utacfC.
Then
\begin{enumerate}[(i)]
\item
	for any $k\le C\root{3}\of{N}$, the number of \krl{s} is
	$$
        O\biggl(\frac{N^2}{k^4}\biggr);
        $$
\item
	more generally, for any $k\le C\root{3}\of{N}$, the number of
	incidences 	between $\P$ and the \krl{s} is $O(N^2/k^3)$;
\item
	the number of incidences between any $M$ straight lines 
	and the $N$ points of $\P$ is 
\[
I= \begin{cases}
	O(M),& \text{if }  N^2<M;\cr
	O(N^{1/2}M^{3/4}),& \text{if } N^{2/3}<M\le N^2;\cr
	O(N^{1/3}M),& \text{if } M\le N^{2/3}.\cr
   \end{cases}
\]
In other words, 
\[
I=O\Bigl(\min\bigl\{M+N^{1/2}M^{3/4}\ ,\ N^{1/3}M\bigr\}\Bigr).
\]
\end{enumerate}
It is also true that these bounds give the best possible 
order of magnitude.
\end{thm}
The forthcoming Section~\ref{ExaSect} provides examples which show that
--- apart from constant factors -- our upper bounds 
cannot be improved.
After some preparatory observations (including our
Main Lemma~\ref{MainLemma}) in Section~\ref{ArrangementSection},
the proof of Theorem~\ref{MainThm} comes in
Section~\ref{ProofSection}.
\end{subsection}

\begin{subsection}{Micha Sharir's ``joints''.}

Given a set of $m$ straight lines $\L=\{L_1,L_2,\ldots,L_m\}$,
a \emph{joint\/} is a point where at least three non-coplanar
$L_i$ meet.
In what follows we denote by $\J(\L)$ the set of joints of $\L$.
\par
It was conjectured by Micha Sharir \cite{Micha:joints}
that if $|\L|=m$ then
$$
|\J(\L)|\le C\cdot m^{3/2},
$$
for an absolute constant $C$.
\par
Here we show the validity of this conjecture if $\J(\L)$ is proper
three--dimen\-sional --- which is true for all known examples with
many joints. Actually, we prove somewhat more (though the original
problem still remains open).
\begin{thm}
If $|\L|=m$ and $\J_0$ is a proper three--dimen\-sional subset
(up to a constant factor $C\ge 1$) of the intersections
--- not necassarily of the joints! --- then
$$
|\J_0|\le C^{3/2} m^{3/2}.
$$
\end{thm}
\proof
Write $n:=|\J_0|$ and consider $C\root3\of n$ planes which cut $\J_0$ into
singletons. Then each $L_i\in\L$ can pass through at most 
$C\root3\of n + 1 \le 2C\root3\of n$ points of $\J_0$, yielding a total of
$m\cdot2C\root3\of n$ incidences.
Since each points in $\J_0$ is incident upon $\ge2$ lines, we have
$2n\le m\cdot2C\root3\of n$, whence
$n\le C^{3/2} m^{3/2}$.
\proofend
\end{subsection}
\end{section}

\begin{section}{Lower bounds}
\label{ExaSect}

\begin{exa}
\label{kRichExample}
For the $N=n^3$ points of an $n\times n\times n$ cube lattice,
and any $2\le k\le n$,
\begin{enumerate}[(a)]
\item
the number of \krl{s} is $\Omega(N^2/k^4)$;
\item
the number of incidences between the lattice points 
and the \krl{s} is $\Omega(N^2/k^3)$.
\end{enumerate}
\end{exa}
[In general, for any $d\ge2$ and  the $N=n^d$ points 
of an $n\times n\times \ldots\times n$ cube lattice
in $\R^d$,  we have at least  $\Omega(N^2/k^{d+1})$
\krl{s} which, of course, produce at least
$\Omega(N^2/k^d)$ incidences.]

\proof
It suffices to show part (a) since it immediately implies part (b).\\
Consider the $N=n^3$ points of $\{1,2,\ldots,n\}^3$.
First we construct  $\Omega(n^3/k^3)$ straight lines 
which all go through the origin $(0,0,0)$, such that
each of them contains approximately $k$ points of the lattice.
These lines will be defined in terms of their points $(u,v,w)$
which is closest to the origin.

We let the coordinates of these points range through
\[
\begin{aligned}
\strut u&= \frac{n}{4k},\ldots,\strut \frac{n}{2k} ;\cr
\strut v&= 1,\ldots,\strut \frac{n}{2k};\cr
\strut w&= 1,\ldots,v\cr
\end{aligned}
\]
such that gcd$(v,w)=1$. 

For each such $(u,v,w)$, the straight line
which passes through it and the origin,
will also pass through at least $2k$ and at most $4k$
points of the cube lattice.
Moreover, the number of such points  $(u,v,w)$ is
\[
\frac{n}{4k} \cdot \sum_{v=1}^{n/(2k)} \phi(v) =
\frac{n}{4k} \cdot \Theta\biggl(\frac{n^2}{k^2}\biggr)
= \Theta\biggl(\frac{n^3}{k^3}\biggr),
\]
where $\phi$ --- i.e., Euler's function ---
gives the number of $w\in\{1\ldots v\}$
which are coprime to $v$, and we used the well-known
fact that $\sum_{i=1}^{m} \phi(i) =\Theta(m^2)$.

Now we shift these lines by each of the vectors 
$(a,b,c) \in \{1,2,\ldots,n/2\}^3$.
Then each new line will still pass through at most $4k$
and, this time, at least $k$ lattice points.
Of course, these $(n/2)^3\cdot \Theta(n^3/k^3)$ lines
are not all distinct. However, each occurs with multiplicity
at most $4k$ whence
\[
\text{number of }k\text{--rich lines }\ge
\biggl(\frac{n}{2}\biggr)^3\cdot\Theta\biggl(\frac{n^3}{k^3}\biggr)
		\cdot \frac{1}{4k} =
\Theta\biggl(\frac{n^6}{k^4}\biggr) =
\Theta\biggl(\frac{N^2}{k^4}\biggr). \proofend
\]
\begin{rem}
\label{HighDimLatticeRem}
A similar construction, with coordinates
$u_1,u_2,\ldots,u_d$ (in place of $u,v,w$),
ranging through
\[
\begin{aligned}
\strut u_1&= \frac{n}{4k},\ldots,\frac{n}{2k} ;\cr
\strut u_2,u_3,\ldots,u_{d-1}&= 1,\ldots,\frac{n}{2k};\cr
\strut u_d&= 1,\ldots,u_{d-1}\cr
\end{aligned}
\]
such that gcd$(u_{d-1},u_d)=1$, gives $\Theta({N^2}/{k^{d+1}})$
lines for a $d$--dimensional $n\times n\times \ldots\times n$
cube lattice with $N=n^d$ points.
\end{rem}
\begin{exa}
To show that the bounds in part (iii) of the 
Main Theorem~\ref{MainThm} are best possible
for all $M$ and $N$, we again consider the
$N=n^3$ points of and $n\times n\times n$
cube lattice.
\begin{enumerate}[(a)]
\item
If $M>N^2/16$, we just draw $M$ lines, each through
at least one point of the lattice.
\item
For $M<N^{2/3}=n^2$, we pick any $M$ of the $n^2$
lattice lines parallel to, say, the $x$--axis.
\item
If $N^2/16>M\ge N^{2/3}$ then we define
\[
2\le k \defeq \frac{N^{1/2}}{M^{1/4}} \le 
	\frac{N^{1/2}}{(N^{2/3})^{1/4}} = N^{1/3} =n
\]
\end{enumerate}
and consider the \krl{s} of the lattice.
According to Example~\ref{kRichExample}.(b),
the number of incidences between these lines 
and the lattice points is 
\[
\Omega\biggl(\frac{N^2}{k^3}\biggr) =
\Omega\biggl(\frac{N^2}{N^{3/2}/M^{3/4}}\biggr) =
\Omega(N^{1/2}M^{3/4}). \proofend
\]
\end{exa}
\end{section}

\begin{section}{Arrangements of planes in $\R^3$.}
\label{ArrangementSection}

\begin{subsection}{Distances and neighborhoods.}

Let $H$ be a (finite) set of $n$ planes in $\R^3$
(or, in general, of hyperplanes in $\R^d$), as in 
Section~\ref{ProperDimDefSection}.
If they are in general position then ---
as it was already mentioned there ---  they cut the space into at most 
$\sim n^3$ open convex cells
(and into $\sim n^d$ in $\R^d$).
The set of these cells, together with their vertices, edges, and faces, 
is called the \emph{arrangement\/} defined by $H$. 
We shall denote it by $\A(H)$.

For two cells $\Cell_i,\Cell_j\in \A(H)$, a natural notion 
of distance is
$$
\dist(\Cell_i,\Cell_j)\defeq\#\{h\in H\ ;\ h \text{ separates }	
		\Cell_i \text{ and } \Cell_j\}.
$$
A spectacular representation is the following: 
pick two points $P_i\in\Cell_i$, $P_j\in\Cell_j$ and connect them
by a straight line segment.
Then the foregoing distance equals the number of $h\in H$ which 
cut the segment $\segm{P_iP_j}$.

It is easy to see that ``dist'' is a metric, 
i.e. it satisfies the triangle inequality.

Our goal is to bound from above --- in terms of $|H|$ ---
the number of pairs $(\Cell_i,\Cell_j)$ whose distance is
at most a given $\rho>0$.
This will be achieved in the Main Lemma~\ref{MainLemma}.

To this end, we define the \rnh of a cell $\Cell_j$ by
$$
B_{\rho}(\Cell_j) = \{\Cell_i\in \A(H)\ ;\ 
		\dist(\Cell_i,\Cell_j)\le \rho	\},
$$
and we note that the number of (ordered) ``$\rho$--close pairs'' 
mentioned above equals 
$$
\sum_{\Cell_j\in\A(H)}	|B_{\rho}(\Cell_j)|.
$$
The next two subsections recall two well-known results,
related to the foregoing \rnh{s} in some sense.
Our main tool (Lemma~\ref{MainLemma}.) comes after these.
\end{subsection}

\begin{subsection}{Zones}

For any (hyper)plane $h\in H$, the \emph{zone\/} of $h$
is the set of cells which ``touch'' $h$, i.e., which 
have a face on $h$.
Also the \emph{\rzn\/}  of $h$ can be defined as the set of
cells $\Cell_j$ for which there is another cell $\Cell_i$ 
in the zone of $h$ for which 
$\dist(\Cell_i,\Cell_j)\le\rho$.
(In this sense the 0--zone coincides with the original zone of $h$.)

\begin{thm}[Matou\v{s}ek]
\label{ZoneThm}
The number of vertices (and, consequently, 
that of the cells, faces, edges)  in the \rzn
of any $h\in H$  is $O(\rho|H|^2)$ in $\R^3$
and  $O(\rho|H|^{d-1})$ in $\R^d$.
\end{thm}
\proof see \cite{Jirka:84} for the bound on the number of vertices.
The rest is implied by the fact that --- according to the 
``general position'' assumption --- each other object
has a vertex furthest from $h$ and each vertex is counted
a bounded number of times (which, of course, depends on the dimension).
\medskip\par
We also re-state this result in terms of \rnh{s}.
To do so, we shall say  that a (hyper)plane $h\in H$
and a \rnh $B_{\rho}(\Cell_j)$ are
\emph{incident upon each other} if $h$ contains 
at least one face of at least one cell 
$\Cell_i\in B_{\rho}(\Cell_j)$.
(It does not matter whether this face is located in the 
``interior'' of the \rnh or on its boundary.)

The next result says that only $O(\rho)$ (hyper)planes
are incident upon an ``average'' neighborhood.
More precisely, we have the following.

\begin{cor}
\label{SumOfnjCor}
For each $\Cell_j\in\A(H)$, denote by $n_j$ the number
of $h\in H$ which are incident upon $B_{\rho}(\Cell_j)$.
Then
$$
\sum_{\Cell_j\in\A(H)} n_j = O(\rho|H|^3)
$$
in $\R^3$ and $O(\rho|H|^d)$ in $\R^d$.
\end{cor}
\proof
Note that $h$ is incident upon $B_{\rho}(\Cell_j)$ 
iff $\Cell_j$ is in the \rzn of $h$.
The rest is just double--counting, using Theorem~\ref{ZoneThm}.
\medskip\par
We also state yet another consequence which can be considered
as the ``younger brother'' (i.e., 2--dimensional version)  of
the forthcoming Main Lem\-ma~\ref{MainLemma}.
\begin{cor}
\label{PlanarBdCor}
In $\R^2$ we have
$$
\sum_{\Cell_j\in\A(H)}	|B_{\rho}(\Cell_j)|
          = O(\rho^2|H|^2).
$$
\end{cor}
\proof
Instead of summing the number of cells $\Cell_i$ in each 
$B_{\rho}(\Cell_j)$, we double--count the triples $(\Cell_j,h,\Cell_i)$ 
such that $h\in H$ bounds  $\Cell_i\in B_{\rho}(\Cell_j)$ and
separates it from $\Cell_j$.\\
On the one hand, the number of these triples cannot be smaller 
than the sum in question (each pair of cells is counted at least once).\\
On the other hand, for a fixed straight line $h\in H$ and a $\Cell_j$ 
in the \rzn
of $h$, the number of the $\Cell_i$ to be counted is at most $2\rho+1$.
(Any two such cells are at distance $\le2\rho$ apart, along the line $h$.)
Thus, using Theorem~\ref{ZoneThm} for $d=2$, we have
$$
\sum_{\Cell_j\in\A(H)}	|B_{\rho}(\Cell_j)|
  \le \# \text{ of triples } \le |H|\cdot O(\rho|H|)\cdot(2\rho+1)
          = O(\rho^2|H|^2).\proofend
$$
As for the second moment $\sum|B_{\rho}(\Cell_j)|^2$, 
it may not always be bounded by a quadratic function of $|H|$
(e.g., if the lines all surround a regular polygon then each of its $|H|$
triangular neighbours has $\ge|H|$ other cells in its $\rho=2$--neighborhood.).
\begin{probl}
\label{RefinementProbl}
Let $\A(H)$ be a simple  arrangement in $\R^2$.
Is it true that it can be refined to an $\A(H^+)$
by adding $O(\sqrt{|H|})$ new straight lines such that
$\sum_{\Cell_j\in\A(H^+)}|B_{\rho}(\Cell_j)|^2= O(\rho^4|H|^2)$?
\end{probl}
It may well be true that one can even force the stronger upper bound
$|B_{\rho}(\Cell_j)|=O(\rho^2)$ for all $\Cell_j\in\A(H^+)$
--- but it is ``even more'' unknown.
\end{subsection}

\begin{subsection}{Levels.}

During this subsection, we study arrangements 
located in a fixed Cartesian coordinate system
and consider the positive half of the $z$--axis
(or that of the $x_d$--axis in $\R^d$) as pointing
``up''.
Thus we can say that a point is ``below'' or ``above'' 
a non-vertical (hyper)plane.

Also, while speaking about levels (to be defined immediately),
we shall assume that none of the (hyper)planes are vertical.

The \emph{level\/}  of a cell $\Cell_j\in\A(H)$ is the number
of $h\in H$ which lie below $\Cell_j$.
This can also be visualized by picking a point $P\in\Cell_j$
and drawing a ray from $P$ downward;
the level of $\Cell_j$ is the number of $h\in H$ which cut this ray.

\begin{thm}[Clarkson]
\label{LevelThm}
The number of vertices, edges, faces, and cells of level $\le\rho$ is
$O(\rho^2|H|)$ in $\R^3$ and
$O(\rho^{\lceil d/2\rceil}|H|^{\lfloor d/2\rfloor})$ in $\R^d$.
\end{thm}
\proof
see \cite{Cl:88} and also Theorem 6.3.1 in \cite{JirkaBook} for vertices;
for the rest proceed as in the proof of  Theorem~\ref{ZoneThm}.
\medskip\par
From now on, we stop stating results for dimensions exceeding three.
The reason for this is that the higher dimensional versions of the
forthcoming bounds --- though usually sharp --- do not seem strong
enough for extending our Main Lemma~\ref{MainLemma} to $d\ge4$.
\begin{cor}
\label{MaxNeighborhoodCor}
In $\R^3$, for any $\Cell_j\in\A(H)$ and $\rho>0$, we have
$$
|B_{\rho}(\Cell_j)|=O(\rho^2|H|).
$$
\end{cor}
\proof
First we pick a point $P\in\Cell_j$ and apply a projective transform 
$\pi$ which maps $P$ to the point at infinity of the $z$--axis.
Consequently, since no $h\in H$ contains $P$, no plane $h$ will be
mapped into vertical position.

For any cell $\Cell_i\in B_{\rho}(\Cell_j)$ and any point
$P_i\in\Cell_i$, the segment $\segm{P_iP}$ intersects $\le\rho$
planes $h\in H$. Moreover, it is mapped to a vertical ray
emanating from $\pi(P_i)$ which, of course, can point either 
downward or upward.

In the former case, the image $\pi(\Cell_i)$ is at level
$\le\rho$ in $\A(\pi(H))$.
According to Theorem~\ref{LevelThm}, there are $O(\rho^2|H|)$
such cells.

Otherwise, in the latter case, we reflect $\pi(H)$ and the
arrangement about the $x$--$y$ plane and apply the same Theorem
to the reflected image.

To sum up, the number of cells in $B_{\rho}(\Cell_j)$
is at most twice the bound in Theorem~\ref{LevelThm}, 
which still makes $O(\rho^2|H|)$. \proofend

\begin{cor}
\label{njSizeCor}
If $B_{\rho}(\Cell_j)$ is incident upon $n_j$ planes $h\in H$
then
$$
|B_{\rho}(\Cell_j)|=O(\rho^2n_j).
$$
\end{cor}
\proof
Those planes which are not incident upon the \rnh
cannot affect its size; we can just delete them
and then apply Corollary~\ref{MaxNeighborhoodCor}.\proofend
\end{subsection}

\begin{subsection}{Graphs of short distances.}

Given an arrangement $\A(H)$ and a $\rho>0$, we define a graph
$G_{\le\rho}$ on the cells $\Cell_j\in\A(H)$ as vertices
(one can visualize them as representative points $P_j\in\Cell_j$)
and edge set $\Er$ by connecting two cells $\Cell_i$, $\Cell_j$
(or, equivalently, the points $P_i$ and $P_j$) by an edge if
$\dist(\Cell_i$, $\Cell_j)\le\rho$.
Our prime tool bounds the number of edges of this graph
in terms of $\rho$ and $|H|$.
\begin{lem}[Main Lemma]
\label{MainLemma}
In $\R^3$, we have
$$
|\Er|=O(\rho^3|H|^3).
$$
\end{lem}
\proof
As in Corollaries~\ref{SumOfnjCor} and \ref{njSizeCor},
denote by $n_j$ the number of planes $h\in H$ 
which are incident upon a cell $\Cell_j\in\A(H)$.
Then
\[
\begin{aligned}
|\Er|&= \sum_{\Cell_j\in\A(H)} |B_{\rho}(\Cell_j)|
		= \sum_{\Cell_j\in\A(H)} O(\rho^2 n_j) = \cr
		&= \rho^2 O\biggl(\sum_{\Cell_j\in\A(H)} n_j\biggr)
		= \rho^2 O(\rho|H|^3) = O(\rho^3|H|^3).\proofend
\end{aligned}
\]
\end{subsection}
\end{section}

\begin{section}{Proof of the Main Theorem~\ref{MainThm}.}
\label{ProofSection}

We demonstrate parts (i)--(iii) one by one, following (and
suitably adapting) an ingenious idea of J.~Solymosi
\cite{SJozsi:???}.
\medskip\par
\proofof{part (i)}
Assume that  a set $\P\subset\R^3$ of $N$ points 
can be cut into singletons by a set $H$ of 
$n\le C\root3\of N$ planes.
In other words, each cell contains at most one point 
$P\in\P$. Moreover, let $2\le k \le n$ be arbitrary.

Define 
\[
\rho=\frac{3n}{k}.
\]
We shall 
make use of the graph $G_{\le\rho}$ of pairs of cells 
$\Cell_i,\Cell_j\in \A(H)$, for which
$\dist(\Cell_i,\Cell_j)\le\rho$.

First we consider a $k$--rich line $l$ and assume that 
the points of $l\cap P$ are $P_1,P_2,\ldots,P_k$
in this linear order.
This $l$ intersects each of the $n$ planes $h\in H$ 
at most once.
Therefore, at most $k/3$ of the segments between consecutive
pairs of points $P_iP_{i+1}$ will intersect more than $\rho$ planes
--- otherwise there would be strictly more than
$(k/3)\cdot(3n/k)=n$ intersections.

Hence there remain at least
\[
(k-1) - \frac{k}{3} = \frac{3k-3-k}{3} \ge \frac{k}{6}
\]
segments which are cut by $\le\rho$ planes $h\in H$.
In terms of the graph $G_{\le\rho}$, each \krl
contributes at least $k/6$ edges.
(Moreover, the latter are all distinct since each
cell contains at most one point and two points 
determine an unique line.)

Since the number of edges satisfies
$|\Er|=O(\rho^3n^3)$ by the Main Lemma~\ref{MainLemma},
we have
\[
\begin{aligned}
\text{number of }k\text{--rich lines } &\le
\frac{O(\rho^3n^3)}{k/6} =
\frac
	{O\biggl(\frac{\ds n^3}{\ds k^3}\cdot n^3\biggr)}
	{\frac{\ds k}{\ds6}}= 
O\biggl(\frac{n^6}{k^4}\biggr)= \cr
&= O\biggl(\frac{N^2}{k^4}\biggr). \proofend \cr
\end{aligned}
\]
\bigskip
\proofof{(ii)}
As a generalization of what was said before, we assume
that a straight line $l$ passes through $k_l\ge k$ points
of the proper 3--dimensional point set $\P$.
Then, just as we have seen, at most $k/3$ (which is 
at most $k_l/3$)  segments will be cut by more than
$\rho=3n/k$ planes $h\in H$, giving way to at least 
\[
(k_l-1) - \frac{k_l}{3} = \frac{3k_l-3-k_l}{3} \ge \frac{k_l}{6}
\]
``close pairs'' and thus at least this many edges of 
$G_{\le\rho}$.
Turning this upside down, for each such line we have that
the number of incidences generated by $l$ is at most six times 
the number of edges of $G_{\le\rho}$ on $l$.
Summing for all \krl{s}, the total number $I$ of incidences 
satisfies
\[
I \le 6\cdot |\Er| = 
6\cdot O\biggl(\frac{n^3}{k^3}\cdot n^3\biggr)= 
O\biggl(\frac{N^2}{k^3}\biggr). \proofend 
\]
\bigskip

\proofof{(iii)}
Consider a set of $N$ points, proper 3--dimensional 
\utacfC.
By definition, the set can be  cut into singletons by
a set $H$ of some $n\le C\root3\of N$ planes.
\medskip

First, for such sets and any $M$ straight lines,
the number of incidences is $O(M\root3\of N)$,
since no line can pass through more than
$C\root3\of N +1 \le (C+1)\root3\of N = O(N^{1/3})$
cells of $\A(H)$, each of which contains at most one
point of the given set.
\medskip

Next, we show another bound which is better than 
the previous one for $M\ge N^{2/3}$.

Denote by $I$ the number of incidences between
our set of $N$ points and $M$ lines.
(Thus an ``average'' line will be incident upon
$\sim I/M$ points.)
\medskip

Put $k=I/(2M)$ and discard all lines which pass 
through less than $k$ points.
Denote by $M\p$ and $I\p$ the number of preserved 
lines  and incidences, respectively.
In total, at most $Mk=I/2$ incidences could be
discarded whence $I\p \ge I/2$.

We distinguish two cases.

Case {\bf I.}\quad
If $k=I/(2M) < 2$ then we have $I< 4M$.
\smallskip

Case {\bf II.}\quad
Otherwise $k=I/(2M) \ge 2$, thus we can apply part (ii)
of the Main Theorem, which yields
\[
I/2 \le I\p = O\biggl(\frac{N^2}{k^3}\biggr) =
	O\biggl(\frac{N^2M^3}{I^3}\biggr),
\]
whence $I^4=O(N^2M^3)$ i.e. $I=O(N^{1/2}M^{3/4})$.
Thus $I=O(M + N^{1/2}M^{3/4})$ anyway, since the
right hand side is an upper bound in either case.
\proofend
\end{section}

\section*{Concluding remarks}

The following questions remain open.

\begin{probl}
Is it true for all $d\ge2$ that if a set $\P\subset\R^d$
of $N$ points 
is proper $d$--dimensional then
\[
\text{number of \krl{s} }= O\biggl(\frac{N^2}{k^{d+1}}\biggr)?
\]
\end{probl}
This order of magnitude, if true, is best possible
(as a function of $N$ and $k$),
as shown by an $N=n\times n \times\ldots\times n$ cube lattice
(see Remark~\ref{HighDimLatticeRem}).
Perhaps a positive answer to the following question could help
in solving the previous problem.
\begin{probl}
Is it true for all $d\ge2$ that  the edge set
$\Er$ of the graph $G_{\le\rho}$ of
``short distances'' defined in terms of 
an arrangement of $n$ hyperplanes  in $\R^d$
satisfies
\[
|\Er|=O(\rho^dn^d)?
\]
\end{probl}
(For $d=1$ the statement is obvious while the cases $d=2$ 
and $d=3$ are Corollary~\ref{PlanarBdCor}
and the Main Lemma~\ref{MainLemma}, respectively.)

\bibliographystyle{alpha}

\end{document}